\documentclass[a4paper,12pt]{article}
\usepackage{amsmath,amsfonts,amssymb,amsthm}

\begin{document}

\newcommand{\R}{\ensuremath{\mathbb{R}}}
\newcommand{\N}{\ensuremath{\mathbb{N}}}
\newcommand{\Z}{\ensuremath{\mathbb{Z}}}

\newcommand{\ep}{\varepsilon}
\newcommand{\lam}{\lambda}
\newcommand{\de}{\delta}
\newcommand{\De}{\Delta}
\newcommand{\al}{\alpha}
\newcommand{\be}{\beta}
\newcommand{\gam}{\gamma}
\newcommand{\Te}{\Theta}

\newcommand{\HH}{{\cal H}}
\newcommand{\II}{{\cal I}}
\newcommand{\NN}{{\cal N}}
\newcommand{\TT}{{\cal T}}
\newcommand{\LL}{{\cal L}}
\newcommand{\DD}{{\cal D}}
\newcommand{\EE}{{\cal E}}
\newcommand{\CCC}{{\cal C}}
\newcommand{\dist}{\mbox{dist}}
\newcommand{\Lip}{\mbox{Lip}}
\newcommand{\Int}{\mbox{Int}}
\newcommand{\NNN}{{\cal N}_n}

{\bf Nonsmooth mappings with Lipschitz shadowing}
\bigskip

Aleksey A. Petrov, Sergei Yu. Pilyugin
\medskip

Faculty of Mathematics and Mechanics, St. Petersburg State University,
University av., 28, 198504, St. Petersburg, Russia



Supported by the Russian Foundation for Basic Research (project 15-01-03797a)
and by St. Petersburg University program ``Stability of dynamical systems
with respect to perturbations and applications to study of applied
problems" (IAS 6.38.223.2014).
\medskip

Abstract. We study conditions under which a piecewise affine mapping has the
Lipschitz shadowing property. As an application, we show that there exists a
homeomorphism with a nonisolated fixed point having the Lipschitz shadowing property.
\bigskip

{\bf 1. Introduction.}
The theory of shadowing of pseudotrajectories (approximate trajectories)
is now a well-developed branch of the theory of dynamical systems
(see, for example, the monographs [1, 2] and the recent survey [3]).

Recently, a lot of attention has been paid to dynamical systems having
special shadowing properties (Lipschitz and H\"older, see [4 -- 6]).
In particular, it was shown in [4] that a diffeomorphism having the
Lipschitz shadowing property is structurally stable (thus, Lipschitz shadowing
property is equivalent to structural stability). The proof in [4]
essentially uses the smoothness of the considered dynamical system
(the Ma\~n\'e theorem [7] giving several characterizations of 
structural stability of diffeomorphisms is applied).

At the same time, it is possible to define the Lipschitz shadowing property
for homeomorphisms (and endomorphisms) of a metric space (see below). 

Of course, if a homeomorphism is topologically conjugate to a structurally
stable diffeomorphism and both the conjugacy and its inverse are uniformly
Lipschitz continuous, then the homeomorphism has the Lipschitz shadowing property.
In this connection, it is
natural to ask: Are homeomorphisms having the Lipschitz shadowing property
similar (in a sense) to structurally stable diffeomorphisms?

In this short note, we give an example of a homeomorphism of the segment
having the Lipschitz shadowing property and a nonisolated fixed point.
This example shows that the answer to the above question is negative.

Let us give the corresponding definitions (for the case of an endomorphism;
for a homeomorphism, the definition is literally the same).

Let $(M,\dist)$ be a metric space and let $f:\;M\to M$ be a continuous 
mapping (we do not distinguish $f$ and the semi-dynamical system generated
by $f$). As usual, a sequence $\pi=\{p_k\in M;\;k\in\Z\}$ is called
a trajectory of $f$ if
$$
p_{k+1}=f(p_k),\quad k\in\Z.
$$

Fix a $d>0$. We say that a sequence $\xi=\{x_k\in M;\;k\in\Z\}$ 
is a $d$-pseudotrajectory of $f$ if
\begin{equation}
\label{0.1}
\dist(x_{k+1},f(x_k))\leq d,\quad k\in\Z.
\end{equation}

The (standard) shadowing property of $f$ means that, given an $\ep>0$, we
can find a $d>0$ such that for any $d$-pseudotrajectory $\xi=\{x_k\}$ of $f$
there is a trajectory $\pi=\{p_k\}$ satisfying the inequalities
\begin{equation}
\label{0.2}
\dist(x_{k},p_k)\leq \ep,\quad k\in\Z.
\end{equation}

Finally, we say that $f$ has the Lipschitz shadowing property if there
exist $\LL,d_0>0$ such that for any $d$-pseudotrajectory $\xi$ of $f$ with $d\leq d_0$
there is a trajectory $\pi$ satisfying inequalities (\ref{0.2}) with
$\ep=\LL d$.
\medskip

The structure of the paper is as follows. In Sec. 2, we prove a general sufficient
condition under which a ``piecewise affine" mapping of $\R^n$ has a ``conditional"
Lipschitz shadowing property (this means that only pseudotrajectories satisfying
some additional assumptions are shadowable). In Sec. 3, we construct the 
above-mentioned example of a homeomorphism of the segment
having the Lipschitz shadowing property and a nonisolated fixed point
(and apply to it the result of Sec. 2). 
\bigskip

{\bf 2. Conditional shadowing result.}
To simplify presentation, we consider a Lipschitz continuous mapping $f:\;\R^n\to \R^n$
with Lipschitz constant $L_0$ (without loss of generality, we assume that $L_0\geq 1$)
for which there exists a family of sets
$G_l\subset \R^n,l\in\Lambda$, with disjoint interiors such that the following conditions hold.

First, for any $l\in\Lambda$ we fix complementary orthogonal 
linear subspaces $S_l$ and $U_l$ of $\R^n$ (let their dimensions be $s_l$ and $u_l$, respectively)
with coordinates $\xi\in S_l$ and $\eta\in U_l$ and denote 
$$
N(\De,p):=\{p+(\xi,\eta):\;|\xi|,|\eta|\leq \De\}
$$
for a point $p\in G_l$ and number $\De>0$.

Let
$$
H_l(\De)=\{p:\;N(\De,p)\subset G_l\}.
$$
\medskip

{\bf Condition~1.} There exists a constant $\lam\in(0,1)$ with the 
following property. For any $l\in\Lambda$ there exist 
$s_l\times s_l$ and $u_l\times u_l$ matrices $A_l$ and $B_l$
 such that
\begin{equation}
\label{0}
\|A_l\|\leq\lam\quad\mbox{and}\quad \|(B_l)^{-1}\|\leq\lam
\end{equation}
and if $p\in H_l(\De)$ for some $\De>0$ (so that $p+(\xi,\eta)\in N(\De,p)$), then 
\begin{equation}
\label{1}
f(p+(\xi,\eta))=f(p)+(A_l\xi,B_l\eta).
\end{equation}
\medskip

{\bf Remark~1. } We impose these simple conditions on the mapping $f$ for the 
following two reasons:

-- they allow us to make the proofs and
estimates maximally ``transparent" (of course, similar results are valid under
more general hyperbolicity conditions on $f$ in the sets $G_l$);

-- precisely these conditions are satisfied in our main application, Theorem~2
below.
\medskip

First we note that the following statement is proved by a standard 
reasoning (for example, it is enough to consider images under the mapping $f$ of the
``rectangles" $N(L_1d,x_{j}),\;0\leq j<m$).
\medskip

{\bf Lemma~1. } {\em Let 
\begin{equation}
\label{l1}
L_1=\frac{1}{1-\lam}.
\end{equation}

If $\{x_k:\;0\leq k\leq m\}$, where $m>0$, is a finite $d$-pseudotrjectory of} $f$ 
({\em this means that inequalities} (\ref{0.1}) {\em are satisfied for} $0\leq k\leq m-1$)
{\em for which there exists an index $l\in\Lambda$ such that
$$
x_{j}\subset H_l(L_1d),\quad 0\leq j<m,
$$
then there exists a point $y$ such that}
\begin{equation}
\label{2}
f^j(y)\in N(L_1 d,x_{j}),\quad 0\leq j\leq m.
\end{equation}

Now we define geometric objects which are important in what follows.
\medskip

Let $p\in G_l,\; l\in\Lambda$; introduce coordinates $(\xi,\eta)$ such that $p$
is the origin and the coordinate subspaces are parallel to $S_l$ and $U_l$, respectively.
Fix $\De_1,\De_2>0$. Consider 
a continuous function $\Xi(\eta)$ that maps the disk
$$
\{\eta:\;\eta\in U_l,\;|\eta|\leq \De_1\}
$$
to $S_l$ and such that
$$
|\Xi(\eta)|\leq \De_2,\quad |\eta|\leq \De_1.
$$
Let $D$ be the graph of $\Xi(\eta)$. 
Denote by ${\cal D}(\De_1,\De_2,p)$ the set of such disks $D$.
\medskip

The following lemma is geometrically obvious.
\medskip

{\bf Lemma~2. } {\em If} $p\in H_l(\De)$, $f(p)\in G_l$, {\em and $D\in
{\cal D}(\De_1,\De_2,p)$, where $\De_1,\De_2\leq \De$, then $f(D)$
contains a disk $D^*$ such that}
$D^*\in{\cal D}(\De_1/\lam,\lam\De_2,f(p))$.
\medskip

{\bf Remark~2. } It is easily seen that in the proof of the main result we
use the statements of Lemmas~1 and 2 (Lipschitz shadowing in $G_l$ with
constant $L_1$ and properties of the images of disks under $f$) plus the ``transversality condition
when we pass from one domain to another" (Condition~2 below).
The assumed linearity of $f$ in the domains $G_l$ just allows us to
make Lemmas~1 and 2 obvious.
\medskip

{\bf Condition~2.} There exist numbers $K\geq L_0+1$ and $\de_0>0$ with the 
following properties. If $L_2=L_1+L_0+1$, $d\leq \de_0$, and there exist three points
$p,q,r$ such that  

(2.1) $p\in G_l$ and $f^2(p)\in G_m$ for some $l,m\in\Lambda$ with $l\neq m$;

(2.2) $q\in H_l(Kd)$ and $r\in H_m(Kd)$;

(2.3) $|p-q|\leq L_1d$ and $|f^2(p)-r|\leq L_2d$;

\noindent 
and

(2.4) $D\in{\cal D}(Kd,d,q)$,

\noindent
then the image $f^2(D)$
contains a disk $D^*$ such that $D^*\in{\cal D}(d,Kd,r)$.
\medskip

{\bf Remark~3. } The above condition is applied in the situation where
points $p$ and $f^2(p)$ belong to different sets $G_l$ and $G_m$ and we
know nothing about the position of the point $f(p)$; in a sense, this condition means
that the image $f^2(D)$ is ``uniformly transverse" to the ``stable subspace"
for $f$ at a point $r$ that is close enough to $f^2(p)$.

Of course, an analog of this condition can be formulated for any pair
of points $p$ and $f^m(p)$, but for our main application (see Sec.~3),
the present form of Condition~2 is enough.
\medskip

We prove the following ``conditional" theorem on Lipschitz shadowing for a mapping
$f$ satisfying the above-formulated conditions. In Theorem~1, we deal with 
finite $d$-pseudotrajectories $\{x_k:\;0\leq k\leq T\}$ of $f$ and show that
there exist $\de_0$ and ${\cal L}$ such that any such finite $d$-pseudotrajectory
with $d\leq\de_0$ 
is ${\cal L}d$-shadowed by a fragment of an exact trajectory of $f$.
It is shown that $\de_0$ and ${\cal L}$ depend only on $f$ and not on the length
of the pseudotrajectory. It is easily seen that if the phase space of a dynamical
system generated by a homeomorphism is locally compact, 
then such a ``finite Lipschitz shadowing property" implies the
Lipschitz shadowing property (cf. [1, Lemma 1.1.1] and the proof of Lemma 4 below).
\medskip

{\bf Theorem~1. } {\em Let }
$X=\{x_k:\;0\leq k\leq T\}$ {\em be a finite $d$-pseudotrajectory of $f$ with $d\leq \de_0$ 
(where $\de_0$ is from Condition} 2). {\em Assume that there
exist (not necessarily different) indices $l_0,l_1,\dots, l_{t}\in\Lambda$
with $l_{i+1}\neq l_i$ and integers
$$
0=m_0<n_0<m_1<n_1<m_2<n_2<\dots<m_t<n_t=T,
$$
where $m_{j+1}=n_j+2,\;j=0,\dots,t-1$, with the following properties:}

(a)
\begin{equation}
\label{5}
x_{k}\in H_{l_j}(K_1d),\quad m_j\leq k\leq n_j,\;j=0,\dots,t, 
\end{equation}
{\em where} $K_1=K+L_1$;

(b)
{\em there exists a positive number $\mu$ for which the inequalities}
\begin{equation}
\label{8}
\mu_j:=n_j-m_{i}\geq\mu,\quad j=0,\dots, t,
\end{equation}

\noindent{\em and}
\begin{equation}
\label{9}
\lam^{\mu}K<1
\end{equation}
{\em are satisfied.}

{\em Let}
\begin{equation}
\label{ll}
\LL=L_0(L_1+2K)+1.
\end{equation}

{\em Then there exists a point $z$ such that}
\begin{equation}
\label{6}
|f^k(z)-x_k|\leq\LL d,\quad k=0,\dots, T.
\end{equation}
\medskip

{\bf Remark~4. } Let us emphasize that only adjacent indices $l_{i+1}$ and $l_i$
are assumed to be different; thus, we do not exclude the situation where
the pseudotrajectory ``returns" to some sets $G_l$ several times.
\medskip

In the proof of Theorem~1, we apply the following statement which is a direct
corollary of Lemma~2.
\medskip

{\bf Lemma~3. } {\em Assume that for some $d>0$ and set $G_l$ there exists a
point $y$ and a number $m$ such that
\begin{equation}
\label{10}
N(Kd,f^k(y))\subset G_l,\quad 0\leq k\leq m,
\end{equation}
and
\begin{equation}
\label{11}
\lam^m K<1.
\end{equation}

Then for any disk $D\in{\cal D}(d,Kd,y)$ there exists a subset $D'\subset D$
such that
\begin{equation}
\label{12}
f^k(D')\subset N(Kd,f^k(y)),\quad 0\leq k\leq m,
\end{equation}
and $f^m(D')$ contains a disk $D^*\in{\cal D}(Kd,d,f^m(y))$}.
\medskip

{\em Proof of Theorem~1. } Fix a $d\leq\de_0$. Condition (a)
allows us to apply Lemma 1 to any ``fragment" 
$$
\{x_{k}:\;m_j\leq k\leq n_j\},\quad j=0,\dots,t,
$$
of the pseudotrajectory $X$ and to find points $y_j,\;j=0,\dots,t,$
such that
\begin{equation}
\label{13}
|f^k(y_j)-x_{m_j+k}|\leq L_1d,\quad 0\leq k\leq \mu_j.
\end{equation}

It follows from condition (\ref{5}) that analogs of inclusions (\ref{10})
in Lemma~3 are satisfied for the points $y_j$:
\begin{equation}
\label{14}
N(Kd,f^k(y_j))\subset G_{l_j},\quad 0\leq k\leq \mu_j,\;j=0,\dots,t.
\end{equation}

Since $\mu_0=n_0-m_0\geq\mu$ (see (\ref{8})), it follows from (\ref{9})
that condition (\ref{11}) of Lemma~3 is satisfied for $y=y_0$ and $m=\mu_0$.

Let $(\xi,\eta)$ be coordinates with coordinate subspaces parallel to
$S_{l_0}$ and $U_{l_0}$, respectively, for which $y_0$ is the origin.

Set 
$$
D_{0,0}=\{(0,\eta):\;|\eta|\leq d\}.
$$
Clearly, $D_{0,0}\in{\cal D}(d,Kd,y_0)$.

Applying Lemma~3, we find a subset $D_0$ of $D_{0,0}$ such that
analogs of inclusions (\ref{12}) are valid, i.e.,
$$
f^k(D_0)\subset N(Kd,f^k(y_0)),\quad 0\leq k\leq\mu_0,
$$
and $f^{\mu_0}(D_0)$ contains a disk $D_0^*\in {\cal D}(Kd,d,f^{\mu_0}(y_0))$.

Let us denote $p=x_{n_0}$, $q=f^{\mu_0}(y_0)$, and $r=y_1$.
It follows from (\ref{13}) (with $j=0$ and $k=n_0$) that
\begin{equation}
\label{16}
|p-q|=|x_{n_0}-f^{\mu_0}(y_0)|=|x_{n_0}-f^{n_0}(y_0)|\leq L_1d.
\end{equation}

Since $X$ is a $d$-pseudotrajectory,
$$
|f^2(p)-x_{m_1}|=|f^2(x_{n_0})-x_{n_0+2}|\leq|f^2(x_{n_0})-f(x_{n_0+1})|+
$$
$$
+|f(x_{n_0+1})-x_{n_0+2}|\leq(L_0+1)d
$$
(recall that $L_0$ is the Lipschitz constant of $f$).
Now we estimate
\begin{equation}
\label{17}
|f^2(p)-r|\leq |f^2(p)-x_{m_1}|+|x_{m_1}-y_1|\leq(L_0+L_1+1)d=L_2d
\end{equation}
(we again refer to (\ref{13}) to estimate the term $|x_{m_1}-y_1|$).

Condition~2 and estimates (\ref{16}) and (\ref{17}) imply that $f^2(D_0^*)$ contains a disk
$D_{1,0}\in{\cal D}(d,Kd,y_1)$. After that, we find a subset $D_{1}\subset D_{1,0}$
that has properties similar to those of $D_0$, and so on.

As a result, we construct sets $D_j,\;j=0,\dots,t$, such that
$$
D_{j+1}\subset f^{\mu_j+2}(D_j),\quad j=0,\dots,t-1,
$$
and 
\begin{equation}
\label{15}
f^k(D_j)\subset N(Kd,f^k(y_j)),\quad 0\leq k\leq\mu_j,\; j=0,\dots,t.
\end{equation}

Hence, there exists a point $z\in D_0$ such that
$$
f^{m_j}(z)\in D_j,\quad j=0,\dots,t.
$$
It follows from inclusions (\ref{15}) and estimates (\ref{13}) that
\begin{equation}
\label{18}
|f^k(z)-x_k|\leq (L_1+2K)d<\LL d,\quad m_j\leq k\leq n_j,\; j=0,\dots,t.
\end{equation}

It remains to estimate the values $|f^k(z)-x_k|$ for $k=n_j+1$.
Let $z'=f^{n_j}(z)$. Then it follows from (\ref{13}) that
$$
|f(z')-x_{n_j+1}|\leq |f(z')-f(x_{n_j})|+|f(x_{n_j})-x_{n_j+1}|\leq
$$
$$
\leq L_0(L_1+2K)d+d=\LL d.
$$
This completes the proof of Theorem~1.
$\Box$
\medskip

{\bf Remark~5. } In parallel to the shadowing property, the so-called
inverse shadowing property is also studied (see, for, example, [8]).
It seems interesting to obtain an analog of Theorem~1 for the Lipschitz
inverse shadowing. Note that the reasoning applied above in the proof
of Theorem~1 cannot be directly transferred to the case of inverse
shadowing.
\bigskip

{\bf 3. Example. }
Consider the segment
$$
I_0=[-7/6,4/3]
$$
and a mapping $f_0:\;I_0\to I_0$ defined as follows:
$$
f_0(x)=1+(x-1)/2,\quad x\in [1/3,4/3],
$$
$$
f_0(x)=2x,\quad x\in (-1/3,1/3).
$$
$$
f_0(x)=-1+(x+1)/2,\quad x\in [-7/6,-1/3].
$$

Clearly, the restriction $f^*$ of $f_0$ to $[-1,1]$ is a homeomorphism of $[-1,1]$
having three fixed points: the points $x=\pm 1$ are attracting and the point $x=0$
is repelling (and this homeomorphism $f^*$ is an example of the so-called
``North Pole -- South Pole" dynamical system; every trajectory starting at a point
$x\neq 0,\pm 1$ tends to an attractive fixed point as time tends to $+\infty$ and 
to the repelling fixed point as time tends to $-\infty$).

Now we define a homeomorphism $f:\;[-1,1]\to [-1,1]$. For an integer $n\geq 0$,
denote $\NNN=2^{-(n+2)}$, and set
\begin{equation}
\label{3.1}
f(x)=\NNN f_0(\NNN^{-1}(x-3\NNN))+3\NNN,\quad x\in(2\NNN,4\NNN].
\end{equation}
This defines $f$ on $(0,1]$. Set $f(0)=0$ and $f(x)=-f(-x)$ for $x\in[-1,0)$.

Clearly, $f$ is a homeomorphism with a nonisolated fixed point $x=0$ 
(for example, every point $x=\pm 2^{-n}$ is
fixed). Let us note that in a neighborhood of any fixed point (with the exception of $x=0$),
$f$ is either linearly expanding with coefficient 2 or linearly contracting with coefficient 1/2.
\medskip

{\bf Theorem 2. } {\em The homeomorphism $f$ has the Lipschitz shadowing property.}
\medskip

Before proving Theorem~2, we prove two auxiliary lemmas 
(and refer to Theorem~1 in the first of them).

In what follows, we denote by $N(d,A)$ the closed
$d$-neighborhood of a set $A$. 
\medskip

{\bf Lemma 4. }{\em The mapping $f_0$ has the Lipschitz shadowing property on $I_0$.}
\medskip

{\em Proof. } 
Let $\xi=\{x_k\in I_0:\;k\in\Z\}$ be a $d$-pseudotrajectory of $f_0$.
In the following (very rough) estimates, we, as usual, decrease values
of $d$, if necessary; every time, the chosen value of $d$ is not more
than the previous values. First we assume that $d\leq d_1<1/24$.

Note that
$$
f_0(-7/6)=-13/12,\quad f_0(4/3)=7/6.
$$

Set
$$
I_0'=\left[-27/24,29/24\right].
$$
Since $\xi$ must belong to $N(d,f_0(I_0))$,
we conclude that
\begin{equation}
\label{3.2}
\xi\subset I_0'.
\end{equation}

Now let us describe the possible position of $\xi$ in $I_0'$.

We note that 
$$
f_0(5/12)=17/24.
$$
If there exists an index $k$ such that $|x_k|\geq 5/12$, then 
$$
|x_{k+i}|>16/24>5/12,\quad i\geq 1.
$$

Note that both $f_0$ and $f_0^{-1}$ have Lipschitz constant 2.
Thus, if $\xi$ is a $d$-pseudotrajectory of $f_0$, then 
$\xi$ is a $2d$-pseudotrajectory of $f_0^{-1}$.

If there exists an index $k_0$ such that $1/4\leq |x_{k_0}|\leq 5/12$, then 
$$
|f_0^{-1}(x_{k_0})|\in[1/8,5/24],
$$
and there exists a $d_2>0$ such that if $d\leq d_2$, then
\begin{equation}
\label{3.3}
|x_{k}|\leq 5/24+2d,\quad k<k_0.
\end{equation}

Thus, only one of the
folowing possibilities can be realized for $d<d_2$:

(1) $|x_k|\leq 1/4$ for $k\in\Z$;

(2) $5/12\leq |x_k|\leq 29/24$ for $k\in\Z$;

(3) there exists an index $k_0$ such that $1/4\leq |x_{k_0}|\leq 5/12$
and inequalities (\ref{3.3}) hold.

In cases (1) and (2), $\xi$ belongs to a domain in which $f_0$ is hyperbolic
(and $\xi$ is uniformly separated from the boundaries of the domain); by Lemma~1, there
exists a $d_3>0$ such that if $d<d_3$, then $\xi$ is $2d$-shadowed by an exact
trajectory of $f_0$.

Consider the remaining case (3) (and let, for definiteness, $k_0=1$ and $x_1>0$; the case 
$x_1<0$ is treated similarly). 

Denote $p=x_0$. Set $G_0=[-1/3,1/3]$ and $G_1=[1/3,29/24]$.
Thus, $p\in G_0$.

As was mentioned, we can take $L_0=2$ and the statement of Lemma~1 holds for $G_0$
and $G_1$ with $L_1=2$.

Since $p\in[1/8-2d,5/24+2d]$, there exists a $d_4>0$ such that if $d<d_4$, then
\begin{equation}
\label{3.4}
5/24+4d<1/3\mbox{ and }N(5d,f^2(p))\subset G_1.
\end{equation}

Take a point $q$ such that
$$
|p-q|\leq L_1d=2d.
$$
In this case, it follows from (\ref{3.4}) that
$q\in G_0$, and, defining disks from $\DD(\De_1,\De_2,q)$,
we must take $S_0=\{0\}$ and $U_0=\R$. Thus, if
$K>0$, then the set $\DD(Kd,d,q)$ contains precisely one disk
$$
D=[q-Kd,q+Kd].
$$
If $K>2$, then $D$ contains the disk
$$
D'=[p-K'd,p+K'd],
$$
where $K'=K-2$.

Clearly, $f^2(D')$ contains the disk
$$
D''=[f^2(p)-K'd/4,f^2(p)+K'd/4].
$$
If a point $r$ satisfies the inequality
$$
|f^2(p)-r|\leq L_2d=5d,
$$
it follows from the second inclusion in (\ref{3.4}) that
$r\in G_1$, and, defining disks from $\DD(\De_1,\De_2,r)$,
we must take $U_1=\{0\}$ and $S_1=\R$. Thus, 
the set $\DD(d,Kd,r)$ consists of points $r'$ such that $|r'-r|\leq d$.

It follows that Condition~2 is satisfied if $d_0\leq d_4$ and
$K'/4\geq 6$.
Thus, it is enough to take $K=26$.

Now, when $L_0$, $L_1$, and $K$ are fixed, it is easily seen 
that there exists a $d_0>0$ such that if $d\leq d_0$, then
\begin{equation}
\label{3.5}
x_k\in H_0(K_1d)=[-1/3+K_1d,1/3-K_1d],\quad k\leq 0,
\end{equation}
and
\begin{equation}
\label{3.6}
x_k\in H_1(K_1d)=[1/3-K_1d,29/24+K_1d],\quad k\geq 2.
\end{equation}

To apply Theorem~1, we fix a natural number $n$ and change indices of
points of the $d$-pseudotrajectory $\xi=\{x_k\}$ to obtain a 
$d$-pseudotrajectory $\xi^{(n)}=\{x_k^{(n)}\}$,
where
$$
x_k^{(n)}=x_{k-n},\quad k\in\Z.
$$

Setting $m_0=0,\;n_0=n,\;m_1=n+2,m_2=2n+2$, we conclude from inclusions
(\ref{3.5}) and (\ref{3.6}) that
$$
x_k^{(n)}\in H_0(K_1d),\quad m_0\leq k\leq n_0,
$$
and
$$
x_k^{(n)}\in H_1(K_1d),\quad m_1\leq k\leq n_1.
$$
Thus, condition (a) of Theorem~1 is satisfied.

It is clear that if $2^{-n-1}K<1$, then condition (b) of Theorem~1 is satisfied
as well.

By Theorem~1, there exists a point $z^{(n)}$ such that
$$
|f^k(z^{(n)})-x_k^{(n)}|\leq\LL d,\quad 0\leq k\leq 2n+2.
$$
Hence, if $\zeta^{(n)}=f^n(z^{(n)})$, then
\begin{equation}
\label{3.7}
|f^k(\zeta^{(n)})-x_k|\leq\LL d,\quad -n\leq k\leq n+2.
\end{equation}

Let $\zeta$ be a limit point of the sequence $\zeta^{(n)}$.
Passing to the limit as $n\to\infty$ in 
(\ref{3.6}) and taking into account that $f$ is a homeomorphism
(so that any $f^k$ is continuous), we see that
$$
|f^k(\zeta)-x_k|\leq\LL d,\quad k\in\Z.
$$
$\Box$
\medskip

The following statement is almost obvious.
\medskip

{\bf Lemma 5. } {\em Let $g$ be a mapping of a segment $J$ and let
numbers $M>0$ and $m$ be given. Consider the mapping
$$
g'(y)=M^{-1}g(M(y-m))+m
$$
on the set
$$
J'=\{y:\;M(y-m)\in J\}.
$$
If $g$ has the Lipschitz shadowing property with constants $\LL,d_0$,
then $g'$ has 
the Lipschitz shadowing property with constants} $\LL,M^{-1}d_0$.
\medskip

{\em Proof. } First we note that if $\{y_k\}$ is a $d$-pseudotrajectory of $g'$
with $d\leq d_0/M$
and $x_k=M(y_k-m)$, then
$$
g(x_k)-x_{k+1}=M(g'(y_k)-y_{k+1}).
$$
Hence, $\{x_k\}$ is an $Md$-pseudotrajectory of $g$.

Since $Md\leq d_0$, there exists a point $p$ such that
$$
|g^k(p)-x_k|\leq \LL Md.
$$

Set $p'=M^{-1}p+m$. Then, obviously,
$$
|(g')^k(p')-y_k|=M^{-1}|g^k(p)-x_k|\leq \LL d.
$$
$\Box$
\medskip

Let us prove Theorem~2.

For a natural $n$, define the segment
$$
I_n=[\al_n,\be_n]=[11\NNN/6,13\NNN/3]
$$
and note that formula (\ref{3.1}) defining $f$ for $x\in(2\NNN,4\NNN]$
is, in fact, valid for $x\in I_n$.

It follows from the equalities
$$
f(\al_n)=23\NNN/12,\quad f(\be_n)=25\NNN/6
$$
that $f(I_n)\subset N(\NNN/12,I_n)$.
Thus, 
if $d<\de(n)=\NNN/12$ and $\{x_k\}$ is a $d$-pseudotrajectory of $f$
that intersects $I_n$, then $\{x_k\}\subset I_n$.

Let $d_0$ and $\LL$ be the constants of Lipschitz shadowing for $f_0$
given by Lemma~4. Since $d_0<1/12$, it follows from Lemma~5 that if 
$\{x_k\}$ is a $d$-pseudotrajectory of $f$
that intersects $I_n$ for some $n>0$, then $\{x_k\}$ is $\LL d$-shadowed.
Of course, a similar statement holds for the segments $I'_n=[-\be_n,-\al_n]$.

To complete the proof, consider a $d$-pseudotrajectory $\xi=\{x_k\}$ of $f$ with $d\leq d_0$
and find the maximal $n_0$ for which $d<\de(n_0)$. Note that then
$$
d\geq \NN_{n_0+1}/12.
$$

If $\xi$ intersects 
one of the segments $I_n$ or $I'_n$ with $n\leq n_0$, then everything is proved.

Otherwise,
$$
|x_k|\leq \al(n_0)=11\NN_{n_0+1}/3\leq 44d,
$$
and $\xi$ is $44d$-shadowed by the rest point $x=0$.
$\Box$
\bigskip

\section*{Funding}

Supported by the Russian Foundation for Basic Research (project 15-01-03797a)
and by St. Petersburg University program ``Stability of dynamical systems
with respect to perturbations and applications to study of applied
problems" (IAS 6.38.223.2014); the first author is also supported by the Chebyshev Laboratory  (Department of Mathematics and Mechanics, St. Petersburg State University)  under RF Government grant 11.G34.31.0026 and by  JSC "Gazprom Neft" and partially supported by  the Dmitry Zimin  Dynasty Foundation.

{\bf References}
\bigskip

1. S. Yu. Pilyugin, {\em Shadowing in dynamical systems},
     Lecture Notes in Mathematics, Springer, \textbf{1706} (1999).

2. K. Palmer,
    {\em Shadowing in dynamical systems. Theory and applications},
     Kluwer (2000).

3. S. Yu. Pilyugin,
Theory of pseudo-orbit shadowing in dynamical systems,
{\em Differential Equations}, {\bf 47} (2011), 1929-1938.

4. S. Yu. Pilyugin and S. B. Tikhomirov, 
Lipschitz shadowing implies structural stability,
{\em Nonlinearity}, {\bf 23} (2010), 2509-2515.

5. S. B. Tikhomirov, H\"older shadowing on finite intervals,
{\em Ergodic Theory Dynam. Systems} (accepted). arXiv:1106.4053v2.

6. A. A. Petrov and S. Yu. Pilyugin, Shadowing near nonhyperbolic
fixed points, {\em Discrete Contin. Dynam. Syst.}, {\bf 34} (2014), 3761-3772.

7. R. Ma\~n\'e, {\em Characterizations of AS diffeomorphisms},
Lecture Notes in Mathematics, Springer, \textbf{597} (1977), 389-394.

8. S. Yu. Pilyugin, Inverse shadowing by continuous methods, 
{\em Discrete Contin. Dynam. Syst.}, {\bf 8} (2002), 29-38.

\end{document}